\documentclass[11pt]{article}

\usepackage{amssymb,latexsym}
\usepackage{amsmath}
\usepackage{cite}
\usepackage{eucal}
\usepackage[latin1]{inputenc}
\usepackage{bbm}

\title{Morita equivalence bimodules for Wick type star products}

\author{\textbf{Nikolai Neumaier\thanks{Nikolai.Neumaier@physik.uni-freiburg.de}},
  \addtocounter{footnote}{2}
  \textbf{Stefan Waldmann\thanks{Stefan.Waldmann@physik.uni-freiburg.de}}
  \\[0.5cm]
  Fakult{\"a}t f{\"u}r Physik\\
  Albert-Ludwigs-Universit{\"a}t Freiburg\\
  Hermann Herder Stra{\ss}e 3\\
  D 79104 Freiburg\\
  Germany
  }

\date{July 2002\\[0.5cm] FR-THEP 2002/09}

\renewcommand{\mathbb}[1]{\mathbbm{#1}}

\newcommand{\im} {{\mathrm i}}
\newcommand{\eu} {{\mathrm e}}

\newcommand{\cc}[1]      {\overline{{#1}}}

\newcommand{\id}         {\mathsf{id}}

\newcommand{\ad}         {\mathrm{ad}}

\newcommand{\End}        {\mathsf{End}}

\newcommand{\Weyl}       {{\mathcal{W}}}
\newcommand{\WL}         {{\mathcal{W}\!\otimes\!\Lambda^{\!\bullet}}}
\newcommand{\WLE}        {{\mathcal{W}\!\otimes\!\Lambda^{\!\bullet}\!\otimes\!\mathcal{E}}}
\newcommand{\WLEnd}      {{\mathcal{W}\!\otimes\!\Lambda^{\!\bullet}\!\otimes\!\mathcal{E}nd(\mathcal{E})}}
\newcommand{\WLO}        {{\mathcal{W}\!\otimes\!\Lambda}}
\newcommand{\WLL}        {{\mathcal{W}\!\otimes\!\Lambda^{\!\bullet}\!\otimes\!\mathcal{L}}}

\newcommand{\nablaE}     {\nabla^E}
\newcommand{\nablaEnd}   {\nabla^{\End(E)}}

\newcommand{\degs}       {\mathrm{deg}_{\mathrm{s}}}
\newcommand{\dega}       {\mathrm{deg}_{\mathrm{a}}}

\newcommand{\HdR}        {\mathrm{H}_{\mathrm{deRham}}}

\newcommand{\Deltafib}   {\Delta_{\mathrm{fib}}}
\newcommand{\starweyl}   {\mathbin{\star_{\scriptscriptstyle\mathrm{Weyl}}}}
\newcommand{\starwick}   {\mathbin{\star_{\scriptscriptstyle\mathrm{Wick}}}}
\newcommand{\starStrwick}{\mathbin{\star'_{\scriptscriptstyle\mathrm{Wick}}}}
\newcommand{\starawick}  {\mathbin{\star_{\cc{\scriptscriptstyle\mathrm{Wick}}}}}
\newcommand{\stark}      {\mathbin{\star_{\kappa}}}

\newcommand{\biwick}     {\mathbin{\bullet_{\scriptscriptstyle\mathrm{Wick}}}}
\newcommand{\biStrwick}  {\mathbin{\bullet'_{\scriptscriptstyle\mathrm{Wick}}}}

\newcommand{\bik}        {\mathbin{\bullet_{\kappa}}}
\newcommand{\Lbiwick}    {\mathbin{\scriptstyle{\blacklozenge}_{\scriptscriptstyle\mathrm{Wick}}}}
\newcommand{\Lbiawick}   {\mathbin{\scriptstyle{\blacklozenge}_{\cc{\scriptscriptstyle\mathrm{Wick}}}}}

\newcommand{\fibweyl}     {\mathbin{\circ_{\scriptscriptstyle\mathrm{Weyl}}}}
\newcommand{\fibwick}     {\mathbin{\circ_{\scriptscriptstyle\mathrm{Wick}}}}
\newcommand{\fibawick}    {\mathbin{\circ_{\cc{\scriptscriptstyle\mathrm{Wick}}}}}
\newcommand{\fibk}        {\mathbin{\circ_{\kappa}}}

\newcommand{\fibbiwick}   {\mathbin{\diamond_{\scriptscriptstyle\mathrm{Wick}}}}
\newcommand{\fibbiawick}  {\mathbin{\diamond_{\cc{\scriptscriptstyle\mathrm{Wick}}}}}

\newcommand{\Dwick}      {\mathcal{D}_{\scriptscriptstyle\mathrm{Wick}}}
\newcommand{\DEwick}      {\mathcal{D}^E_{\scriptscriptstyle\mathrm{Wick}}}
\newcommand{\DStrwick}      {\mathcal{D'}_{\scriptscriptstyle\mathrm{Wick}}}
\newcommand{\Dawick}     {\mathcal{D}_{\cc{\scriptscriptstyle\mathrm{Wick}}}}

\newcommand{\tauwick}    {\tau_{\scriptscriptstyle\mathrm{Wick}}}
\newcommand{\tauStrwick} {\tau'_{\scriptscriptstyle\mathrm{Wick}}}
\newcommand{\tauEwick}    {\tau^E_{\scriptscriptstyle\mathrm{Wick}}}
\newcommand{\tauawick}   {\tau_{\cc{\scriptscriptstyle\mathrm{Wick}}}}
\newcommand{\rwick}      {r_{\scriptscriptstyle\mathrm{Wick}}}
\newcommand{\rStrwick}    {r'_{\scriptscriptstyle\mathrm{Wick}}}
\newcommand{\rEwick}     {r^E_{\scriptscriptstyle\mathrm{Wick}}}
\newcommand{\rawick}     {r_{\cc{\scriptscriptstyle\mathrm{Wick}}}}

\newtheorem{lemma}{Lemma}
\newtheorem{proposition}{Proposition}
\newtheorem{theorem}{Theorem}
\newtheorem{corollary}{Corollary}

\newtheorem{remark}{Remark}

\newenvironment{proof}{\small{\sc Proof:}}{{\hspace*{\fill} $\square$\\}}

\allowdisplaybreaks

\begin{document}

\maketitle

\begin{abstract}
    In this paper, the notion of star products with separation of
    variables on a Kähler manifold is extended to bimodule
    deformations of (anti-) holomorphic vector bundles over a Kähler
    manifold. Here the Fedosov construction is appropriately adapted
    using the geometric data of a connection in the vector bundle. 
    Moreover, the relation between the star products of Wick
    and anti Wick type is clarified by constructing a canonical Morita
    equivalence bimodule as bimodule deformation of the canonical line
    bundle over the Kähler manifold.
\end{abstract}



%
%

\section{Introduction}
\label{sec:intro}

Deformation quantization as introduced in \cite{bayen.et.al:1978} is a
well-established and successful way of understanding the transition
from classical physics to quantum physics as a deformation of the
algebraic structures, see
\cite{dito.sternheimer:2002,gutt:2000,halbout:2002,waldmann:2002a} for
recent reviews. The classical observable algebra is modelled by the
complex-valued smooth functions $C^\infty(M)$ on the phase space $M$,
which is a Poisson manifold. In particular, $C^\infty(M)$ is a Poisson
algebra. The deformation is done by means of a star product $\star$
which is an associative $\mathbb{C}[[\lambda]]$-bilinear
multiplication for $C^\infty(M)[[\lambda]]$ such that in zeroth order
of the deformation parameter $\lambda$ it coincides with the pointwise
product and in first order of $\lambda$ the $\star$-commutator equals
$\im$ times the Poisson bracket. In case of convergence $\lambda$
corresponds to Planck's constant $\hbar$. Existence and classification
of such products is now well-established
\cite{dewilde.lecomte:1983b,fedosov:1994a,omori.maeda.yoshioka:1991,kontsevich:1997:pre,nest.tsygan:1995a,bertelson.cahen.gutt:1997,weinstein.xu:1998,gutt.rawnsley:1999}.

In this article we shall consider a particular situation, namely when
$M$ is a Kähler manifold and the Poisson structure is induced by its
Kähler symplectic form. Deformation quantization of Kähler manifolds
has a long history dating back to the work of Berezin
\cite{berezin:1975a,berezin:1975b}. The relation between
Berezin-Toeplitz quantization and star products has been explored by
various authors 
\cite{cahen.gutt.rawnsley:1995,cahen.gutt.rawnsley:1994,cahen.gutt.rawnsley:1993,cahen.gutt.rawnsley:1990,schlichenmaier:1997,karabegov.schlichenmaier:2001a,karabegov.schlichenmaier:2001b,landsman:1998}.
Since on a Kähler manifold one has a compatible complex structure it
is meaningful to speak of bidifferential operators (and hence of star
products) which differentiate one function in holomorphic directions
and the other in anti-holomorphic directions only. This property of a
star product was known for various examples and led Karabegov to the
notion of star products with \emph{separation of variables}. In
\cite{karabegov:1996,karabegov:1998a} he proved existence and gave a
classification of such star products. In an alternative approach,
Fedosov's construction was adapted to the Kähler situation and used to
show the existence of such star products as well
\cite{bordemann.waldmann:1997a}, see also \cite{neumaier:2002:pre} for
a classification in this context and \cite{karabegov:2000} for a
comparison of both approaches.

It turns out that on a Kähler manifold one has (at least) three
canonically given star products: the Weyl ordered Fedosov star
product $\starweyl$ as well as the Wick and anti Wick star product
$\starwick$ and $\starawick$, obtained by the above mentioned
modified Fedosov construction, all using the Kähler connection as
starting point. Moreover, the characteristic classes of
these three star products are explicitly known to be
\begin{equation}
    \label{eq:classeswwaw}
    c(\starweyl) = \frac{1}{\im\lambda} [\omega],
    \quad
    c(\starwick) 
    = \frac{1}{\im\lambda} [\omega] - \im \pi c_1(L_{\textrm{can}})
    \quad
    c(\starawick) 
    = \frac{1}{\im\lambda} [\omega] + \im \pi c_1(L_{\textrm{can}}),
\end{equation}
where $L_{\textrm{can}} = \bigwedge^{(n,0)} T^*M$ is the canonical
line bundle of $M$, see
e.g.~\cite{karabegov.schlichenmaier:2001a,hawkins:2000,neumaier:2002:pre}.

In \cite{bursztyn.waldmann:2000b} the notion of deformation
quantization of vector bundles was introduced: If $E \to M$ is a
(complex) vector bundle and $\star$ is a star product for $M$ then a
deformation quantization of $E$ is a deformed right module structure
$\bullet$ for $\Gamma^\infty(E)[[\lambda]]$ with respect to
$\star$. It was shown that such a $\bullet$ always exists and is
unique up to equivalence. Moreover, $\bullet$ also induces an
associative formal deformation $\star'$ of
$\Gamma^\infty(\End(E))[[\lambda]]$ 
together with a corresponding left module structure $\bullet'$ such
that $\Gamma^\infty(E)$ is deformed into a $\star'$-$\star$ bimodule
via $\bullet'$ and $\bullet$. Finally, the bimodule
gives a \emph{Morita equivalence bimodule} for the deformed algebras
$(\Gamma^\infty(\End(E))[[\lambda]], \star')$ and
$(C^\infty(M)[[\lambda]], \star)$. In case of a line bundle $E=L$ one
obtains a star product $\star'$ for 
$C^\infty(M) = \Gamma^\infty(\End(L))$ and finally arrives at the
classification of star products up to Morita equivalence: In the
symplectic case two star products $\star'$ and $\star$ are Morita
equivalent if and only if there is a symplectic diffeomorphism $\psi$
of $M$ such that 
$\psi^*c(\star') - c(\star) \in 2\pi\im \HdR(M,\mathbb{Z})$ is an
\emph{integral} de Rham class, see \cite{bursztyn.waldmann:2002a}.

The aim of this work is two-fold: On one hand we shall construct
particular bimodule deformations $\bullet'$, $\bullet$ for a complex
vector bundle $E$ over $M$ using the ideas of Fedosov's construction
as in \cite{waldmann:2002b} in order to obtain separation of
variables for $\bullet'$ and $\bullet$ as well, if $E$ is (anti-)
holomorphic. On the other hand we give a Fedosov construction of the
deformed bimodule structure of $L_{\textrm{can}}$ which yields a
Morita equivalence bimodule for $\starwick$ and
$\starawick$. According to (\ref{eq:classeswwaw}) and the general
classification theorem in \cite{bursztyn.waldmann:2002a} such
deformations necessarily exist. So our main emphasize is the
canonical and constructive way how it can be obtained.

The paper is organized as follows: In Section~\ref{sec:fedosov} we
collect some basic and well-known results on the Fedosov construction
and adapt them to the Kähler situation. In order to handle the
Weyl ordered, the Wick ordered and the anti Wick ordered case
simultaneously we introduce a one-parameter family of fibrewise
$\kappa$-ordered products. In Section~\ref{sec:bundle} we state the main
theorems of the Fedosov construction for the star products and the
corresponding bimodule multiplications. Up to here the results are
essentially standard. In Section~\ref{sec:wicktype} we investigate the
case of (anti-) holomorphic vector bundles and show the separation of
variables properties of the Wick and anti Wick ordered
products. Section~\ref{sec:local} contains local expressions for the
deformed multiplications which can also be used to characterize them
globally. In Section~\ref{sec:Hermmetrics} we give a deformed version of
a Hermitian fibre metric and investigate its compatibility with the
holomorphic structure. Here we also give local expressions. Finally,
Section~\ref{sec:morita} contains the construction of the canonical
Morita equivalence bimodule structure on $L_{\textrm{can}}$ for
$\starwick$ and $\starawick$. In an appendix we have collected some
standard results on Kähler geometry in order to explain our notation.

\smallskip

\noindent
\textbf{Conventions:} By $C^\infty(M)$ we denote the
\emph{complex-valued} smooth functions and similarly
$\Gamma^\infty(T^*M)$ stands for the complex-valued smooth one-forms
etc. Moreover, we use Einstein's summation convention in local
expressions. Finally, we do not need the positive definiteness of the
Kähler metric. So all results are still valid on 
\emph{semi Kähler manifolds}. However, for ease of notation we shall
not emphasize this in the text.

\medskip
\noindent
\textbf{Acknowledgements:} We would like to thank Martin
Schlichenmaier for valuable discussions on topics in Kähler geometry.

%
%

\section{The Fedosov construction on Kähler manifolds}
\label{sec:fedosov}

In this section we shall briefly recall the set-up for the Fedosov
construction in order to explain our notation where we mainly follow
\cite{waldmann:2002b}. Details and proofs can be found in
Fedosov's book \cite{fedosov:1996}. For the additional structures on a
Kähler manifold we refer to
\cite{bordemann.waldmann:1997a,neumaier:2002:pre} as well as to
Appendix~\ref{sec:kaehler}.

Given a complex vector bundle $E \to M$ over a Kähler manifold $M$ we
define the following $\mathbb{C}[[\lambda]]$-modules
\begin{equation}
    \label{eq:WDef}
    \Weyl := 
    \prod_{s=0}^\infty \Gamma^\infty (\mbox{$\bigvee$}^s T^*M)
    [[\lambda]],
\end{equation}
\begin{equation}
    \label{eq:WLDef}
    \WL 
    := \prod_{s=0}^\infty \Gamma^\infty 
    (\mbox{$\bigvee$}^s T^*M \otimes 
    \mbox{$\bigwedge$}^{\!\bullet}\, T^*M)[[\lambda]],
\end{equation}
\begin{equation}
    \label{eq:WLEDef}
    \WLE := \prod_{s=0}^\infty \Gamma^\infty
    (\mbox{$\bigvee$}^s T^*M \otimes 
    \mbox{$\bigwedge$}^{\!\bullet}\, T^*M \otimes E)[[\lambda]],
\end{equation}
\begin{equation}
    \label{eq:WLEndDef}
    \WLEnd
    := \prod_{s=0}^\infty \Gamma^\infty 
    (\mbox{$\bigvee$}^s T^*M 
    \otimes \mbox{$\bigwedge$}^{\!\bullet}\, T^*M 
    \otimes \End(E))[[\lambda]].
\end{equation}
Then $\WL$ becomes a super-commutative associative algebra where the
fibrewise product is defined by 
$(f\otimes\alpha)(g\otimes\beta)= f\vee g\otimes\alpha\wedge\beta$. In
particular $\Weyl \subseteq \WL$ is a commutative sub-algebra. Using
the fibrewise composition of endomorphisms of $E$ we observe that
$\WLEnd$ becomes an associative algebra as well, which is no longer
super-commutative unless $E=L$ is a line bundle. We can view $\WL$ as
sub-algebra of $\WLEnd$. Finally, $\WLE$ is a bimodule for $\WLEnd$
from the left and for $\WL$ from the right. Besides the usual
symmetric and anti-symmetric form degree $\degs$ and $\dega$, we have
the $\lambda$-degree and the \emph{total degree}, which is twice the
$\lambda$-degree plus the symmetric degree. One also has the operators
$\delta = dx^i \wedge i_s(\partial_i)$ and
$\delta^* = dx^i \vee i_a(\partial_i)$ which satisfy
$\delta^2 = 0 = (\delta^*)^2$. If one defines
$\delta^{-1} a = \frac{1}{k+l} \, \delta^*a$ for homogeneous $a$ with
symmetric degree $k$ and anti-symmetric degree $\l$ such that 
$k+l \ne 0$ and $\delta^{-1}a = 0$ otherwise, then one has the
\emph{Hodge-de Rham decomposition}
\begin{equation}
    \label{eq:hodge}
    \delta \delta^{-1} + \delta^{-1} \delta + \sigma = \id,
\end{equation}
where $\sigma$ denotes the projection onto the part with symmetric and
anti-symmetric degree $0$.

In a next step one needs a symplectic connection $\nabla$ on $M$ and a
connection $\nabla^E$ for $E$. For $\nabla$ we shall always use the
Kähler connection while $\nabla^E$ shall be specified later.
One has the induced connection 
$\nabla^{\End(E)} = [\nabla^E, \cdot]$ on $\End(E)$ and the
connections extend to super-derivations of anti-symmetric degree $+1$
\begin{equation}
    \label{eq:DDef}
    D: \WL \to \WL^{+1},
\end{equation}
\begin{equation}
    \label{eq:DEDef}
    D^E: \WLE \to \WL^{+1}\!\otimes\!\mathcal{E},
\end{equation}
\begin{equation}
    \label{eq:DEndDef}
    D': \WLEnd \to \WL^{+1}\!\otimes\!\mathcal{E}nd(\mathcal{E}).
\end{equation}
Then $D^E$ is a module derivation along $D$ and $D'$, respectively. A
simple computation shows that $\delta$ super-commutes with $D'$, $D^E$
and $D$.

One observes that the curvature $R^E$ of $\nabla^E$ is an element in
$\WLO^2\!\otimes\!\mathcal{E}nd(\mathcal{E})$ with $\degs R^E =
0$. For $R^E$ and for the symplectic curvature tensor $R \in \WLO^2$,
see~(\ref{eq:symplRDef}), we have the Bianchi identities
$\delta R = 0 = \delta R^E$ and $DR = 0 = D'R^E$.

Now we pass to the deformed fibrewise products and bimodule
structures. Originally, Fedosov used the fibrewise Weyl product,
but on a Kähler manifold one also has fibrewise analogues of the Wick
and anti Wick product, see
\cite{bordemann.waldmann:1997a,neumaier:2002:pre,karabegov:2000}.
First one defines the fibrewise operator
\begin{equation}
    \label{eq:FibrePoisson}
    \mathcal{P} 
    = \Lambda^{ij} \, i_s(\partial_i) \otimes i_s(\partial_j)
\end{equation}
acting on $\Weyl\otimes\Weyl$, where 
$\Lambda = \frac{1}{2}\Lambda^{ij}\partial_i\wedge\partial_j$ is the
Poisson tensor in local coordinates. Clearly $\mathcal{P}$ is globally
defined. As we have a complex structure the operators
\begin{equation}
    \label{eq:PccPDef}
    P = g^{k\cc{\ell}} \, i_s(Z_k) \otimes i_s(\cc{Z}_\ell)
    \quad
    \textrm{and}
    \quad
    \cc{P} = g^{k\cc{\ell}} \, i_s(\cc{Z}_\ell) \otimes i_s(Z_k)
\end{equation}
are globally well-defined as well and one has 
$\mathcal{P} = \frac{2}{\im}(P - \cc{P})$. Finally, we need the
fibrewise Laplace operator
\begin{equation}
    \label{eq:DeltafibDef}
    \Deltafib = g^{k\cc{\ell}} \, i_s(Z_k) i_s(\cc{Z}_\ell),
\end{equation}
which again is globally well-defined and satisfies the relation
\begin{equation}
    \label{eq:DeltafibMu}
    \Deltafib \circ \mu 
    = \mu \circ 
    (\Deltafib \otimes \id + P + \cc{P} + \id \otimes \Deltafib),
\end{equation}
where $\mu: \Weyl \otimes \Weyl \to \Weyl$ denotes the fibrewise
undeformed product of $\Weyl$. We extend all the operators
$\mathcal{P}$, $P$, $\cc{P}$, and $\Deltafib$ to $\WL$, $\WLE$, and
$\WLEnd$. Now the fibrewise Weyl product is defined by
\begin{equation}
    \label{eq:fibreWeyl}
    a \fibweyl b 
    := \mu \circ \eu^{\frac{\im\lambda}{2}\mathcal{P}} 
    a \otimes b
\end{equation}
and clearly gives an associative deformation of $\mu$. For
$\kappa \in \mathbb{R}$ we define the 
\emph{fibrewise equivalence transformations}
\begin{equation}
    \label{eq:SkappaDef}
    S = \eu^{\lambda\Deltafib}
    \quad
    \textrm{and}
    \quad
    S^\kappa = \eu^{\lambda\kappa\Deltafib}
\end{equation}
and the fibrewise products
\begin{equation}
    \label{eq:circkappaDef}
    a \fibk b := S^\kappa (S^{-\kappa}a \fibweyl S^{-\kappa} b)
    = \mu \circ \eu^{(\kappa+1)\lambda P + (\kappa-1)\lambda \cc{P}} 
    a \otimes b.
\end{equation}
Clearly $\fibk$ is again a fibrewise associative deformation of $\mu$
and fibrewisely equivalent to $\fibweyl$ via $S^\kappa$. Besides the
Weyl ordered case ($\kappa = 0$) the Wick ordered case ($\kappa = 1$)
and the anti Wick ordered case ($\kappa = -1$) are of particular
interest as here $\cc{P}$ respectively $P$ are absent in
(\ref{eq:circkappaDef}). We shall denote these fibrewise products by
$\fibwick$ and $\fibawick$.

We also observe that the bimodule structure of $\WLE$ can be deformed
yielding a bimodule structure with respect to the deformed products
$\fibk$ of $\WLEnd$ and $\WL$. We shall denote the corresponding
bimodule multiplications by $\fibk$ as well.
\begin{remark}
    \label{remark:almostcomplex}
    It is clear that $\Deltafib$, $P$, $\cc{P}$, and $\fibk$ can be
    defined as soon as one has an almost complex structure compatible
    with the symplectic structure on an arbitrary symplectic manifold.
    In principle one can carry through Fedosov's construction in this
    case as well \cite{karabegov.schlichenmaier:2001b}. However, the
    resulting star products $\star_\kappa$ seem to have no particular
    properties like `separation of variables' for the cases $\kappa =
    \pm 1$. This only will happen on a (semi-) Kähler manifold.
\end{remark}

We finally mention some further relations between the various
operators. First note that $[\Deltafib, \delta] = 0$ whence 
$S^\kappa \delta S^{-\kappa} = \delta$ for all $\kappa$. Moreover,
since $\nabla$ is the Kähler connection we have 
$[\Deltafib, D] = 0$ as well as $[\Deltafib, D^E] = 0$ and
$[\Deltafib, D'] = 0$. This implies $S^\kappa D S^{-\kappa} = D$ and
hence $D$ is a $\fibk$-derivation as well. Analogously, $D'$ and $D^E$
are (module-) derivations with respect to $\fibk$. The following lemma
is a slight variation of \cite[Prop.~4.1]{bordemann.waldmann:1997a}
and \cite[Sect.~5.3]{fedosov:1996}:
\begin{lemma}
    \label{lemma:Dsquare}
    We have 
    $S^\kappa R = R + \kappa\lambda\Deltafib R 
    = R + \kappa\lambda\varrho$ 
    and hence 
    \begin{equation}
        \label{eq:Dsquare}
        D^2 = \frac{\im}{\lambda} \ad_\kappa (S^\kappa R)
        = \frac{\im}{\lambda} \ad_\kappa (R),
    \end{equation}
    as the Ricci form $\varrho$, see (\ref{eq:Ricciform}), is central
    with respect to $\fibk$. Moreover, $S^\kappa R^E = R^E$ and hence 
    \begin{equation}
        \label{eq:DEDERE}
        (D^E)^2 = \frac{\im}{\lambda} \ad_\kappa (R) + R^E
        \quad
        \textrm{and}
        \quad
        (D')^2 = \frac{\im}{\lambda} \ad_\kappa (R - \im\lambda R^E).
    \end{equation}
\end{lemma}

%
%

\section{Fedosov star products and deformed vector bundles}
\label{sec:bundle}

Using the results of the previous section as well as the standard
arguments of Fedosov's construction
\cite{fedosov:1994a,waldmann:2002b} we easily arrive at the following
theorem:
\begin{theorem}
    \label{theorem:FedosovI}
    We fix $\kappa \in \mathbb{R}$. For any series of closed two-forms
    $\Omega_\kappa \in 
    \lambda \Gamma^\infty(\bigwedge^2T^*M)[[\lambda]]$ 
    there exist unique $r_\kappa \in \WLO^1$, 
    $r'_\kappa \in \WLO^1\!\otimes\!\mathcal{E}nd(\mathcal{E})$ of
    total degree $\ge 3$ such that
    \begin{equation}
        \label{eq:rkappa}
        \delta r_\kappa = R + Dr_\kappa + \frac{\im}{\lambda}
        r_\kappa \fibk r_\kappa + \Omega_\kappa
        \quad
        \textrm{and}
        \quad
        \delta^{-1} r_\kappa = 0
    \end{equation}
    and
    \begin{equation}
        \label{eq:rkappaII}
        \delta r'_\kappa = R - \im\lambda R^E + D'r'_\kappa +
        \frac{\im}{\lambda} r'_\kappa \fibk r'_\kappa + \Omega_\kappa
        \quad
        \textrm{and}
        \quad
        \delta^{-1} r'_\kappa = 0.
    \end{equation}
    In this case the super derivations
    \begin{equation}
        \label{eq:DDEnd}
        \mathcal{D}_\kappa 
        = -\delta + D + \frac{\im}{\lambda} \ad_\kappa(r_\kappa)
        \quad
        \textrm{and}
        \quad
        \mathcal{D}'_\kappa 
        = -\delta + D' + \frac{\im}{\lambda} \ad_\kappa(r'_\kappa)
    \end{equation}
    have square zero. The maps
    \begin{equation}
        \label{eq:sigmaI}
        \sigma: \ker \mathcal{D}_\kappa \cap \WLO^0 
        \to C^\infty(M)[[\lambda]]
        \quad
        \textrm{and}
        \quad
        \sigma: \ker \mathcal{D}'_\kappa \cap 
        \WLO^0\!\otimes\!\mathcal{E}nd(\mathcal{E})
        \to \Gamma^\infty(\End(E))[[\lambda]]
    \end{equation}
    are $\mathbb{C}[[\lambda]]$-linear bijections with inverses
    denoted by $\tau_\kappa$ and $\tau'_\kappa$,
    respectively. Finally,
    \begin{equation}
        \label{eq:stark}
        f \stark g = \sigma(\tau_\kappa(f) \fibk \tau_\kappa(g))
        \quad
        \textrm{and}
        \quad
        A \stark' B 
        = \sigma(\tau'_\kappa (A) \fibk \tau'_\kappa (B))
    \end{equation}
    define associative deformations of $C^\infty(M)$ and
    $\Gamma^\infty(\End(E))$, respectively, and $\stark$ is a star
    product with characteristic class
    \begin{equation}
        \label{eq:classstark}
        c(\stark) = \frac{1}{\im\lambda} 
        ([\omega] + [\Omega_\kappa] -\kappa\lambda[\varrho]).
    \end{equation}
\end{theorem}
\begin{proof}
    The proof is completely standard and follows from
    \cite{fedosov:1994a,waldmann:2002b} with some obvious
    modifications. We only indicate the computation of the
    characteristic class: We define $r = S^{-\kappa} r_\kappa$ whence
    clearly
    \begin{equation}
        \label{eq:deltar}
        \delta r = R - \kappa\lambda\varrho + Dr + \frac{\im}{\lambda}
        r \fibweyl r + \Omega_\kappa
    \end{equation}
    using Lemma~\ref{lemma:Dsquare} and $\Deltafib\Omega_\kappa = 0$.
    Moreover, $\delta^{-1}r$ is some element of total degree $\ge 3$.
    Thus the characteristic class of the Weyl ordered Fedosov star
    product $\star$ built out of $r$ is simply given by the expression
    in (\ref{eq:classstark}), see e.g.~\cite{neumaier:1999:pre}. But
    $\star$ and $\stark$ are equivalent as their Fedosov derivatives
    are fibrewisely conjugate by $S^\kappa$,
    see~\cite[Prop.~1]{bordemann.neumaier.waldmann:1998}.
\end{proof}

For the deformed bimodule structure of $\Gamma^\infty(E)$ with respect
to $\stark'$ and $\stark$ we proceed completely analogously to
\cite{waldmann:2002b}. First we define
\begin{equation}
    \label{eq:rEkappa}
    r^E_\kappa := \frac{\im}{\lambda}(r'_\kappa - r_\kappa),
\end{equation}
which is a well-defined formal power series in $\lambda$ as
$r'_\kappa$ and $r_\kappa$ coincide in zeroth order of
$\lambda$. Then we define 
$\mathcal{D}^E_\kappa: \WLE \to \WL^{+1}\!\otimes\!\mathcal{E}$ by
\begin{equation}
    \label{eq:DEkappa}
    \mathcal{D}^E_\kappa = -\delta + D^E +
    \frac{\im}{\lambda}\ad_\kappa(r_\kappa) + r^E_\kappa,
\end{equation}
where $\ad_\kappa(r_\kappa)$ is defined as usual and $r^E_\kappa$ acts
by $\fibk$-left multiplications. The following construction is an
immediate adaption of \cite[Thm.~3]{waldmann:2002b} to our
$\kappa$-ordered situation:
\begin{theorem}
    \label{theorem:FedosovII}
    The Fedosov derivative $\mathcal{D}^E_\kappa$ satisfies
    \begin{equation}
        \label{eq:BiModDerI}
        \mathcal{D}^E_\kappa (\Psi \fibk b) 
        = \mathcal{D}^E_\kappa \Psi \fibk b
        + (-1)^k \Psi \fibk \mathcal{D}_\kappa b,
    \end{equation}
    \begin{equation}
        \label{eq:BiModDerII}
        \mathcal{D}^E_\kappa (a \fibk \Psi) 
        = \mathcal{D}'_\kappa a \fibk \Psi +
        (-1)^\ell a \fibk \mathcal{D}^E_\kappa \Psi,
    \end{equation}
   and $(\mathcal{D}^E_\kappa)^2 = 0$, for 
   $a \in \WLO^\ell\!\otimes\!\mathcal{E}nd(\mathcal{E})$,
   $\Psi \in \WLO^k\!\otimes\!\mathcal{E}$, and 
   $b \in \WL$. Moreover, 
   $\sigma: \ker\mathcal{D}^E_\kappa \cap 
   \WLO^0\!\otimes\!\mathcal{E} \to \Gamma^\infty(E)[[\lambda]]$
   is a $\mathbb{C}[[\lambda]]$-linear bijection with inverse denoted
   by $\tau^E_\kappa$. Finally,
   \begin{equation}
       \label{eq:bulletkappaLeft}
       A \bullet'_\kappa s := \sigma(\tau'_\kappa(A) \fibk \tau^E_\kappa(s))
   \end{equation}
   \begin{equation}
       \label{eq:bulletkapparight}
       s \bik f := \sigma(\tau^E_\kappa(s) \fibk \tau_\kappa(f))       
   \end{equation}
   defines a $\stark'$-$\stark$ bimodule deformation of
   $\Gamma^\infty(E)$.
\end{theorem}
\begin{remark}
    \begin{enumerate}
    \item In the particular case of a line bundle $E=L$ we obtain a
        star product $\star'_\kappa$ for 
        $C^\infty(M)[[\lambda]] = \Gamma^\infty(\End(L))[[\lambda]]$
        with a corresponding bimodule deformation. Then the
        characteristic class of $\star'_\kappa$ is given by
        \begin{equation}
            \label{eq:classstarleft}
            c(\star'_\kappa) = \frac{1}{\im\lambda}([\omega] +
            [\Omega_\kappa] - \kappa\lambda[\varrho]) 
            + 2\pi\im c_1(L),
        \end{equation}
        where $c_1(L)$ is the Chern class of $L$. This follows either
        from \cite[Thm.~3.1]{bursztyn.waldmann:2002a} or by an
        analogous argument as in \cite[Cor.~2]{waldmann:2002b}.
    \item In general, $\star'_\kappa$ and $\stark$ are Morita
        equivalent deformations and 
        $(\Gamma^\infty(E)[[\lambda]],\bullet'_\kappa, \bullet_\kappa)$
        is a Morita equivalence bimodule
        \cite[Prop.~1]{waldmann:2002b}.
    \end{enumerate}
\end{remark}

%
%

\section{The Wick Type Properties of $\starwick$, $\starStrwick$,
$\biwick$ and $\biStrwick$}
\label{sec:wicktype}

In this section we shall consider the deformations $\star_\kappa$,
$\star'_\kappa$ and the bimodule structures $\bullet_\kappa$,
$\bullet'_\kappa$ more closely in the case $\kappa=1$. Actually an
analogous consideration can also be carried out in the case $\kappa =
-1$ by almost trivial generalizations of the given results.  We set
$\Omega = \Omega_1 = \Omega_{\scriptscriptstyle\textrm{Wick}}$.

First we need some additional notations that make use of the complex
structure that enable us to consider the splittings into holomorphic
and anti-holomorphic part of the mappings involved in the Fedosov
construction. For a detailed discussion of this topic the reader is
referred to \cite[Appx.~A]{neumaier:2002:pre}. By $\pi_z$ we denote the
projection onto the holomorphic form part in the symmetric and the
anti-symmetric part of $\WL$. Then $\pi_z$ naturally extends to a
projection defined on $\WLEnd$ and $\WLE$ as well. Analogously
$\pi_{\cc{z}}$ denotes the projection onto the anti-holomorphic form
part of $\WL$ that also extends naturally to $\WLEnd$ and $\WLE$.
Obviously we have $\sigma = \pi_z \pi_{\cc{z}}=\pi_{\cc{z}} \pi_z$.
From the very definitions of the products and the fibrewise bimodule
multiplications $\fibwick$ it is easy to see that
\begin{equation}
    \label{eq:pizcirccomp}
    \pi_z(F \fibwick G) = \pi_z((\pi_z F)\fibwick G)
    \quad\textrm{and}\quad
    \pi_{\cc{z}}(F \fibwick G) =
    \pi_{\cc{z}}(F\fibwick (\pi_{\cc{z}}G)),
\end{equation}
where $F,G$ are elements in $\WL$, $\WLE$ or $\WLEnd$ such that the
muliplications make sense. Using a local holomorphic chart of $M$
it is easy to see that $\delta =
\delta_z +
\delta_{\cc{z}}$, $\delta^*= {\delta_z}^{\!\!\!*}+
{\delta_{\cc{z}}}^{\!\!\!*}$, $D = D_z +D_{\cc{z}}$, $D^E = D^E_z +
D^E_{\cc{z}}$ and $D'= D'_z + D'_{\cc{z}}$, where for instance
$D'_z(f \otimes \alpha\otimes A):=
\nabla_{Z_i}f \otimes dz^i\wedge\alpha \otimes A+ f \otimes
\partial \alpha\otimes A + f \otimes dz^i\wedge\alpha\otimes
\nablaEnd_{Z_i}A$ and $D'_{\cc{z}}(f \otimes
\alpha\otimes A):=\nabla_{\cc{Z}_i}f \otimes d\cc{z}^i\wedge\alpha
\otimes A+ f \otimes \cc{\partial} \alpha\otimes A + f \otimes
d\cc{z}^i\wedge\alpha\otimes\nablaEnd_{\cc{Z}_i}A$. The other
splittings are defined similarly. Completely analogously to the
definition of $\delta^{-1}$ one defines ${\delta_z}^{\!\!\!-1}$ for
$a$ with $dz^i\vee i_s(Z_i)a=k a$ and $ dz^i
\wedge i_a(Z_i)a = l a$ by ${\delta_z}^{\!\!\!-1}a :=
\frac{1}{k+l} {\delta_z}^{\!\!\!*}a$ in case $k+l \neq 0$ and
${\delta_z}^{\!\!\!-1}a := 0$ in case $k+l=0$.
${\delta_{\cc{z}}}^{\!\!\!-1}$ is defined in the analogous way and
an easy computation yields the following decompositions:
\begin{equation}
    \label{eq:HodgeholaholZer}
    {\delta_z}^{\!\!\!-1} \delta_z  
    + \delta_z {\delta_z}^{\!\!\!-1} +
    \pi_{\cc{z}}=\id
    \quad\textrm{and}\quad
    {\delta_{\cc{z}}}^{\!\!\!-1} \delta_{\cc{z}}  
    + \delta_{\cc{z}} {\delta_{\cc{z}}}^{\!\!\!-1}  +
    \pi_z =\id.
\end{equation}
Furthermore we have the relations: $\pi_z \delta = \delta_z \pi_z$,
$\pi_z \delta^{-1}={\delta_z}^{\!\!\!-1}\pi_z$, $\pi_z D = D_z
\pi_z$, $\pi_z D^E = D^E_z \pi_z$ and $\pi_z D' = D'_z \pi_z$ and
analogous formulas with $\cc{z}$ instead of $z$. Using the equation
$\delta^2=0$, the fact that $\delta$ super-commutes with $D$, $D^E$
and $D'$ and the equations for $D^2$, ${D^E}^2$ and ${D'}^2$ it is
easy to derive super-commutation relations for the holomorphic and
the anti-holomorphic parts of the involved mappings, see also 
\cite[Lem.~3]{neumaier:2002:pre}.
\begin{lemma}
    Let $\rwick$, $\rStrwick$ and $\rEwick$ denote the elements of 
    $\WLO^1$, $\WLO^1\!\otimes \mathcal{E}nd(\mathcal E)$ constructed
    in the preceding section in case $\kappa =1$. Then we have the
    following:
    \begin{enumerate}
    \item
        In case $\Omega$ is of type $(1,1)$ then in addition
        \begin{equation}\label{eq:rWiProjNull}
            \pi_z \rwick =0 
            \quad\textrm{and}\quad
            \pi_{\cc{z}}\rwick=0.
        \end{equation}
    \item
        In case $\Omega$ and $R^E$ are of type $(1,1)$ then
        \begin{equation}\label{eq:rWiStrProjNull}
            \pi_z \rStrwick=\pi_{\cc{z}}\rStrwick=0 
            \quad\textrm{and}\quad
            \pi_z \rEwick = \pi_{\cc{z}} \rEwick=0.
        \end{equation}
    \end{enumerate}
\end{lemma}
\begin{proof}
The proof of i.) can be found in
\cite[Lem.~4.1]{neumaier:2002:pre}. Applying 
$\pi_z$ to (\ref{eq:rkappaII}) we obtain
\[
\delta_z \pi_z \rStrwick = D'_z \pi_z \rStrwick + \frac{\im}{\lambda}
\pi_z ((\pi_z \rStrwick)\fibwick \rStrwick)\qquad \textrm{ and } \qquad
{\delta_z}^{\!\!\!-1} \pi_z \rStrwick = 0
\]
using $\pi_z R^E=\pi_z\Omega=0$. Using the decomposition
${\delta_z}^{\!\!\!-1} \delta_z +
\delta_z {\delta_z}^{\!\!\!-1} + \pi_{\cc{z}}=\id$ this implies that
$\pi_z \rStrwick$ is a fixed point of the mapping $L :
\pi_z(\WLO^1\! \otimes\mathcal{E}nd(\mathcal E)) \ni a \mapsto
{\delta_z}^{\!\!\!-1}( D'_z a + \frac{\im}{\nu}\pi_z(a\fibwick
\rStrwick))\in \pi_z(\WLO^1\!\otimes\mathcal{E}nd(\mathcal E))$ that raises the
total degree at least by one and hence has a unique fixed point.
Obviously $L$ has $0$ as trivial fixed point implying $\pi_z
\rStrwick=0$ by uniqueness. Analogously one proves that
$\pi_{\cc{z}}\rStrwick=0$. Then the statement for $\rEwick$
follows.
\end{proof}

From the explicit shape of the fibrewise products and the bimodule
multiplications $\fibwick$ it is obvious that the products
$\starwick$, $\starStrwick$ and the bimodule muliplications
$\biwick$, $\biStrwick$ are completely determined by the knowledge
of the projections on the totally holomorphic resp.
anti-holomorphic part of the respective Fedosov-Taylor series.
Therefore it is worth deriving simpler formulas for these
projections in order to investigate the products and the bimodule
multiplications.
\begin{proposition}
    \label{proposition:TaySerProjEq}
    Let $R^E$ and $\Omega$ be of type $(1,1)$ then the projections of
    $\tauwick(f)$, $\tauStrwick(A)$ and $\tauEwick(s)$ for $f\in
    C^\infty(M)[[\lambda]]$, $A\in
    \Gamma^\infty(\mathsf{End}(E))[[\lambda]]$ and $s \in
    \Gamma^\infty(E)[[\lambda]]$ satisfy the equations:
    \begin{eqnarray}\label{eq:ztauWi}
        \delta_z \pi_z \tauwick(f) = D_z \pi_z \tauwick(f) -
        \frac{\im}{\lambda} \pi_z((\pi_z \tauwick(f))\fibwick \rwick)
        &\textrm{and}& \sigma(\pi_z \tauwick(f)) =f,\\
        \label{eq:ccztauWi}
        \delta_{\cc{z}} \pi_{\cc{z}}\tauwick(f) = D_{\cc{z}} \pi_{\cc{z}}
        \tauwick(f) +\frac{\im}{\lambda}\pi_{\cc{z}}(\rwick \fibwick
        (\pi_{\cc{z}} \tauwick(f)))&\textrm{and}&
        \sigma(\pi_{\cc{z}} \tauwick(f)) =f,\\
        \label{eq:ztauStrWi}
        \delta_z \pi_z \tauStrwick(A) = D'_z \pi_z \tauStrwick(A) -
        \frac{\im}{\lambda} \pi_z ((\pi_z\tauStrwick(A))\fibwick \rStrwick)
        &\textrm{and}&\sigma(\pi_z\tauStrwick(A))=A,\\
        \label{eq:ccztauStrWi}
        \delta_{\cc{z}}\pi_{\cc{z}} \tauStrwick(A) = D'_{\cc{z}} \pi_{\cc{z}}
        \tauStrwick(A) + \frac{\im}{\lambda}\pi_{\cc{z}}(\rStrwick \fibwick
        (\pi_{\cc{z}}\tauStrwick(A)))&\textrm{and}&
        \sigma(\pi_{\cc{z}}\tauStrwick(A))=A,\\
        \label{eq:ztauEWi}
        \delta_z \pi_z \tauEwick(s) = D^E_z \pi_z \tauEwick
        (s) -\frac{\im}{\lambda} \pi_z((\pi_z
        \tauEwick(s))\fibwick \rwick)
        &\textrm{and}&\sigma(\pi_z \tauEwick(s))= s,\\
        \label{eq:ccztauEWi}
        \delta_{\cc{z}}\pi_{\cc{z}}\tauEwick(s) = D^E_{\cc{z}}
        \pi_{\cc{z}}
        \tauEwick (s) + \frac{\im}{\lambda} \pi_{\cc{z}}(\rStrwick
        \fibwick
        (\pi_{\cc{z}}\tauEwick(s))) &\textrm{and}&
        \sigma(\pi_{\cc{z}}\tauEwick(s))= s
    \end{eqnarray}
    from which they are uniquely determined.
\end{proposition}
\begin{proof}
    The proof of equations (\ref{eq:ztauWi}) and (\ref{eq:ccztauWi})
    again can be found in \cite[Prop.~4.3]{neumaier:2002:pre}. It is easy
    to see that this proof can be slightly modified to obtain the
    remaining statements of the proposition.
\end{proof}

As an easy consequence of the preceding proposition we find:
\begin{lemma}
    \label{lemma:projundef}
    Let $\mathcal{U} \subseteq M$ be an open subset of M and let $R^E$
    and $\Omega$ be of type $(1,1)$.
    \begin{enumerate}
    \item For all $f\in C^\infty(M)$ anti-holomorphic on $\mathcal{U}$
        we have
        $\pi_z \tauwick(f)\big|_{\mathcal{U}} = f\big|_{\mathcal{U}}$.
    \item For all $g\in C^\infty(M)$ holomorphic on $\mathcal{U}$ we
        have $\pi_{\cc{z}} \tauwick(g)\big|_{\mathcal{U}} =
        g\big|_{\mathcal{U}}$.
    \item For all $A\in \Gamma^\infty(\End(E))$ such that
        $\nablaEnd_{Y} A\big|_{\mathcal{U}} =0$ for all
        $Y\in\Gamma^\infty(TM^{1,0})$ we have $\pi_z
        \tauStrwick(A)\big|_{\mathcal{U}} = A\big|_U$.
    \item For all $B\in \Gamma^\infty(\End(E))$ such that $\nablaEnd_X
        B\big|_{\mathcal{U}} = 0$ for all
        $X\in\Gamma^\infty(TM^{0,1})$ we have $\pi_{\cc{z}}
        \tauStrwick(B)\big|_{\mathcal{U}} = B\big|_{\mathcal{U}}$.
    \item For all $s\in \Gamma^\infty(E)$ such that $\nablaE_{Y}
        s\big|_{\mathcal{U}}=0$ for all $Y\in\Gamma^\infty(TM^{1,0})$
        we have $\pi_z \tauEwick(s)\big|_{\mathcal{U}} =
        s\big|_{\mathcal{U}}$.
    \item For all $t\in \Gamma^\infty(E)$ such that $\nablaE_X
        t\big|_{\mathcal{U}}=0$ for all $X\in\Gamma^\infty(TM^{0,1})$
        we have $\pi_{\cc{z}}\tauEwick(t)\big|_{\mathcal{U}} =
        t\big|_{\mathcal{U}}$.
    \end{enumerate}
\end{lemma}
\begin{proof}
    To prove these assertions one just has to show that the given
    expressions for the projections of the respective Fedosov-Taylor
    series solve the equations given in Proposition
    \ref{proposition:TaySerProjEq} which is an easy task using that
    the projections of $\rwick$ and $\rStrwick$ onto the
    totally holomorphic and anti-holomorphic form part vanish.
\end{proof}

From the results of the preceding lemma we can deduce that the
products $\starwick$, $\starStrwick$ and the bimodule
multiplications $\biwick$, $\biStrwick$ have the following Wick
type properties:
\begin{theorem}\label{theorem:WickType}
    Let $\mathcal{U} \subseteq M$ be an open subset of M and let $R^E$ and
    $\Omega$ be of type $(1,1)$.
    \begin{enumerate}
    \item For all $s\in \Gamma^\infty(E)$, $f, g \in C^\infty(M)$ such
        that $g$ is holomorphic in $\mathcal{U}$ we have
        \begin{equation}\label{eq:holbulstarvonrechts}
            s\biwick g \big|_{\mathcal{U}} = s g\big|_{\mathcal{U}}
            \quad\textrm{and}\quad 
            f \starwick g \big|_{\mathcal{U}} = f g\big|_{\mathcal{U}}.
        \end{equation}
    \item For all $s \in \Gamma^\infty(E)$, 
        $A, B \in \Gamma^\infty(\End(E))$ such that 
        $\nablaEnd_{Y} B\big|_{\mathcal{U}} =0$ for all 
        $Y \in \Gamma^\infty(TM^{1,0})$ we have 
        \begin{equation}
            \label{eq:antiholbulStrstarStrvonlinks}
            B \biStrwick s\big|_{\mathcal{U}} = Bs\big|_{\mathcal{U}}
            \quad\textrm{and}\quad 
            B \starStrwick A \big|_{\mathcal{U}} = B A\big|_{\mathcal{U}}.
        \end{equation}
    \item For all $t, s\in \Gamma^\infty(E)$, 
        $A\in\Gamma^\infty(\End(E))$,
        $f \in C^\infty(M)$ such that $\nablaE_X t\big|_{\mathcal{U}} =0$ for all 
        $X \in \Gamma^\infty(TM^{0,1})$ and $\nablaE_{Y} s\big|_{\mathcal{U}} =0$ for
        all $Y \in \Gamma^\infty(TM^{1,0})$ we have
        \begin{equation}\label{eq:holbulStrvrechtsantiholbulvlinks}
            A \biStrwick t\big|_{\mathcal{U}} = A t\big|_{\mathcal{U}}
            \quad\textrm{and}\quad 
            s \biwick f \big|_{\mathcal{U}} = s f\big|_{\mathcal{U}}.
        \end{equation}
    \item For all $A, B \in \Gamma^\infty(\End(E))$, 
        $f, g \in C^\infty(M)$
        such that $\nablaEnd_X B\big|_{\mathcal{U}}=0$ for all 
        $X \in \Gamma^\infty(TM^{0,1})$ and 
        $g$ anti-holomorphic on $\mathcal{U}$ we have
        \begin{equation}\label{eq:holstarStrvrechtsaholstarvlinks}
            A \starStrwick B \big|_{\mathcal{U}} = A B\big|_{\mathcal{U}}
            \quad\textrm{and}\quad 
            g\starwick f \big|_{\mathcal{U}} = g f\big|_{\mathcal{U}}.
        \end{equation}
    \end{enumerate}
\end{theorem}
To conclude this section we want to discuss some more concrete
situations where the precondition that the curvature endomorphism
$R^E$ is of type $(1,1)$ naturally occurs. If $(E,h)$ is a (anti-)
holomorphic vector bundle with Hermitian fibre metric then there
exists a unique connection $\nablaE$ which is compatible with the
(anti-) holomorphic structure as well as with the Hermitian fibre
metric (cf. Appendix~\ref{sec:kaehler}). In this case it is known that
$R^E$ is of type $(1,1)$. Moreover, in this case the conditions of the
preceding Proposition under which the products and the bimodule
multiplications are just the pointwise ones just mean that the
respective sections resp. functions are locally holomorphic resp.
anti-holomorphic. In the case of a holomorphic vector bundle $(E,h)$ a
section $s$ is called locally anti-holomorphic on $\mathcal{U}$ if the
section $s^\flat$ of the dual bundle $E^*$ which is also a holomorphic
vector bundle defined by $s^\flat(s'):= h(s, s')$ is locally
holomorphic on $\mathcal{U}$. Similarly a section $A \in
\Gamma^\infty(\End(E))$ is called locally anti-holomorphic if the
section $A^*$ defined by $h(A^*s,s')=h(s, A s')$ for all $s, s' \in
\Gamma^\infty(E)$ is locally holomorphic. In the case of an
anti-holomorphic vector bundle $(E,h)$ one proceeds similarly to
define the notion of locally holomorphic sections.

%
%

\section{Local Expressions for the Bimodule Multiplications
$\biwick$, $\biStrwick$}
\label{sec:local}

In this section we want to consider the bimodule multiplications
constructed above in the case of a (anti-) holomorphic vector
bundle $E$ of fibre dimension $k$ equipped with a connection that
is compatible with the (anti-) holomorphic structure and the
curvature of which is of type $(1,1)$. We denote by
$\{\mathcal{U}_\alpha\}$ a good open cover of $M$.
\begin{proposition}
    In case $E$ is an anti-holomorphic vector bundle and $e_\alpha$ is
    an anti-holomorphic frame on $\mathcal U_\alpha$ one has
    \begin{equation}\label{eq:piztauEWiloc}
        \pi_z \tauEwick(s)\big|_{\mathcal U_\alpha} =
        e_\alpha \pi_z \tauwick(s_\alpha)= \pi_z \tauwick(
        s^i_\alpha)\otimes e_{\alpha,i},
    \end{equation}
    where $s\in \Gamma^\infty(E)[[\lambda]]$ has been written as
    $s|_{\mathcal U_\alpha}= e_\alpha s_\alpha =
    e_{\alpha,i}s^i_\alpha$. Hence we have for all
    $f\in C^\infty(M)[[\lambda]]$ 
    \begin{equation}\label{eq:bulWiloc}
        s \biwick f \big|_{\mathcal U_\alpha} = e_\alpha (s_\alpha
        \starwick f)=  e_{\alpha,i}(s^i_\alpha \starwick f).
    \end{equation}
    Conversely, $\biwick$ is globally well-defined by equation
    (\ref{eq:bulWiloc}) and thus it is completely determined by
    $\starwick$.
\end{proposition}
\begin{proof}
    To prove that $\pi_z \tauEwick(s)$ is locally given by $\pi_z
    \tauwick(s^i_\alpha)\otimes e_{\alpha,i}$ it is enough to show
    that this expression solves the equations (\ref{eq:ztauEWi}) what
    is easily done using that the frame $e_\alpha$ is
    anti-holomorphic.  But then (\ref{eq:bulWiloc}) is immediate from
    the definition of $\biwick$. It remains to show that $\biwick$ is
    globally defined by this equation. To this end let $e_\beta =
    e_\alpha \phi_{\alpha\beta}$ another anti-holomorphic frame where
    $\phi_{\alpha\beta}$ denotes the anti-holomorphic transition
    matrix on $\mathcal U_\alpha\cap \mathcal U_\beta$ as in
    Appendix~\ref{sec:kaehler}. We find
    \[
    e_\alpha (s_\alpha \starwick f) = e_\beta \phi_{\beta \alpha}
    ((\phi_{\alpha\beta}s_\beta)\starwick f) = e_\beta (s_\beta
    \starwick f)
    \]
    on $\mathcal U_\alpha\cap\mathcal U_\beta$, where we have used
    that the star product $\starwick$ is of Wick type implying
    $(\phi_{\alpha\beta}s_\beta)\starwick f=
    \phi_{\alpha\beta}(s_\beta\starwick f)$ since the entries of
    $\phi_{\alpha\beta}$ are anti-holomorphic.
\end{proof}
\begin{remark}
    In fact the preceding proposition states that for an
    anti-holomorphic vector bundle $E$ with a connection that is
    compatible with the anti-holomorphic structure the right module
    multiplication $\biwick$ is canonical in the sense that is
    independent of the connection. Actually the statement of the
    proposition could also be proved without using the local formula
    for $\pi_z \tauEwick(s)$ but only the Wick type properties of
    $\biwick$ and $\starwick$ according to Theorem
    \ref{theorem:WickType}. To do so observe that
    $\Gamma^\infty(E)[[\lambda]] \ni s = e_\alpha s_\alpha = e_\alpha
    \biwick s_\alpha$ since $e_\alpha$ is an anti-holomorphic frame.
    But then
    \[
    s \biwick f|_{\mathcal U_\alpha} = (e_\alpha \biwick s_\alpha)
    \biwick f = e_\alpha \biwick (s_\alpha \starwick f ) =
    e_\alpha (s_\alpha
    \starwick f).
    \]
    To see that this yields a global definition for $\biwick$ we note
    the following equations that are valid on $\mathcal U_\alpha\cap
    \mathcal U_\beta$
    \begin{equation}
        \label{eq:niceframes}
        e_\beta = e_\alpha \phi_{\alpha\beta} = e_\alpha
        \biwick \phi_{\alpha \beta}
        \quad\textrm{and}\quad
        s_\alpha = 
        \phi_{\alpha \beta} s_\beta 
        = \phi_{\alpha \beta} \starwick s_\beta 
    \end{equation}
    \begin{equation}
        \label{eq:nicematrix}
        \phi_{\alpha\beta}
        \phi_{\beta\alpha}= \phi_{\alpha\beta}\starwick
        \phi_{\beta\alpha}= \id 
        \quad\textrm{and}\quad
        \phi_{\alpha\beta}
        \phi_{\beta\gamma} \phi_{\gamma\alpha}=\phi_{\alpha\beta}\starwick
        \phi_{\beta\gamma} \starwick \phi_{\gamma\alpha}=\id,
    \end{equation}
    where $\phi_{\alpha \beta}$ denotes the \emph{classical} transition
    matrices, which imply the statement using
    \cite[Lem.~9.3]{waldmann:2002a}. 
\end{remark}

As an immediate consequence of these observations we have:
\begin{corollary}
    Let $E$ be an anti-holomorphic vector bundle and let
    ${\mathbin{\tilde\bullet_{\scriptscriptstyle\mathrm{Wick}}}}$ be a
    $\starwick$ right module multiplication on
    $\Gamma^\infty(E)[[\lambda]]$ that has the Wick type property,
    then ${\mathbin{\tilde\bullet_{\scriptscriptstyle\mathrm{Wick}}}}$
    coincides with $\biwick$.
\end{corollary}

In contrast the bimodule multiplication $\biStrwick$ drastically
simplifies in case $E$ is a holomorphic vector bundle, but even in
this case it actually depends on the connection $\nablaE$ and hence is
not canonical like $\biwick$. Nevertheless it is completely determined
by the product $\starStrwick$ as we state it in the following
proposition.
\begin{proposition}
    In case $E$ is a holomorphic vector bundle and $e_\alpha$ is a
    holomorphic frame and $e^\alpha$ is the dual frame of the vector
    bundle $E^*$ on $\mathcal U_\alpha$ one has
    \begin{equation}
        \label{eq:piccztauEWiloc}
        \pi_{\cc{z}} \tauEwick (s)\big|_{\mathcal U_\alpha} =
        \frac{1}{k}
        \pi_{\cc{z}} \tauStrwick(s \otimes e^{\alpha,i}) e_{\alpha,i}
    \end{equation}
    for all $s \in \Gamma^\infty(E)[[\lambda]]$. Hence we have for all
    $A\in \Gamma^\infty(\End(E))[[\lambda]]$
    \begin{equation}
        \label{eq:bulStrWiloc}
        A \biStrwick s\big|_{\mathcal U_\alpha} =
        \frac{1}{k} (A \starStrwick ( s
        \otimes e^{\alpha,i})) e_{\alpha,i}.
    \end{equation}
    Moreover, $\biStrwick$ is globally well-defined by equation
    (\ref{eq:bulStrWiloc}) and thus it is completely determined by
    $\starStrwick$.
\end{proposition}
\begin{proof}
    An easy computation using that $e^\alpha$ is a holomorphic frame
    of $E^*$ shows that the expression given in
    (\ref{eq:piccztauEWiloc}) solves the equations
    (\ref{eq:ccztauEWi}). But then the local formula
    (\ref{eq:bulStrWiloc}) follows immediately. In fact this could
    also be proved directly just using the Wick type property of
    $\starStrwick$ and $\biStrwick$. To see this we observe that on
    $\mathcal U_\alpha$ we have
    \[
    s = \frac{1}{k} (s\otimes e^{\alpha,i}) e_{\alpha,i}
    = \frac{1}{k} (s\otimes e^{\alpha,i})
    \biStrwick e_{\alpha,i}
    \]
    since $e_\alpha$ is a holomorphic frame. But this implies that on
    $\mathcal U_\alpha$
    \begin{eqnarray*}
        A \biStrwick s &=& \frac{1}{k} A
        \biStrwick ((s\otimes e^{\alpha,i}) \biStrwick e_{\alpha,i})=
        \frac{1}{k} (A\starStrwick(
        s\otimes e^{\alpha,i})) \biStrwick e_{\alpha,i}\\ &=&
        \frac{1}{k} (A\starStrwick( s\otimes
        e^{\alpha,i}))e_{\alpha,i}.
    \end{eqnarray*}
    We now write $e_{\alpha,i}= {\phi_{\alpha\beta}}^j_i e_{\beta,j}$
    and $e^{\alpha,i}={\phi^{\alpha\beta}}^i_k e^{\beta,k}$ with
    ${\phi^{\alpha\beta}}^n_r {\phi_{\alpha\beta}}^r_i= \delta^n_i$.
    Using the Wick type property of $\starStrwick$ and $\biStrwick$
    and that the transition matrices $\phi_{\alpha\beta},
    \phi^{\alpha\beta}$ have holomorphic entries we moreover get
    \begin{eqnarray*}
        (A\starStrwick( s\otimes e^{\alpha,i}))e_{\alpha,i} &=&
        (A\starStrwick( s\otimes e^{\beta,k})\starStrwick
        ({\phi^{\alpha\beta}}^i_k\id))\biStrwick ({\phi_{\alpha\beta}}^j_i
        e_{\beta,j})\\ &=& A \biStrwick ( (( s\otimes e^{\beta,k})
        \biStrwick (({\phi^{\alpha\beta}}^i_k\id)\biStrwick
        ({\phi_{\alpha\beta}}^j_i e_{\beta,j})))) = (A\starStrwick(
        s\otimes e^{\beta,j}))e_{\beta,j}
    \end{eqnarray*}
    on $\mathcal U_\alpha\cap \mathcal U_\beta$ and hence $\biStrwick$
    is globally well-defined by (\ref{eq:bulStrWiloc}).
\end{proof}

\begin{corollary}
    Let $E$ be a holomorphic vector bundle and let
    ${\mathbin{\tilde\bullet'_{\scriptscriptstyle\mathrm{Wick}}}}$ be
    a $\starStrwick$ left module multiplication on
    $\Gamma^\infty(E)[[\lambda]]$ that has the Wick type property,
    then
    ${\mathbin{\tilde\bullet'_{\scriptscriptstyle\mathrm{Wick}}}}$
    coincides with $\biStrwick$.
\end{corollary}

%
%

\section{Deformed Hermitian Metrics}
\label{sec:Hermmetrics}

In this section we consider a Hermitian fibre metric $h$ for $E$ which
is assumed to be either a holomorphic or anti-holomorphic vector bundle.
Moreover, let $\nablaE$ be the canonical connection for $(E,h)$.
As in the preceding sections we assume $\Omega$ is of type $(1,1)$
and in addition we consider the case where $\Omega$ is real
i.e. $\cc{\Omega} = \Omega$. 
Now $h$ equips $\Gamma^\infty(\End(E))$ with a natural $^*$-involution
defined by 
$h(A  s, s') = h(s, A^* s')$ for $A\in\Gamma^\infty(\End(E))$ and 
$s, s' \in \Gamma^\infty (E)$. 
We extend this involution, together with the complex conjugation, to a
super-$^*$-involution of $\WLEnd$ and $\WL$, respectively. It is
well-known that the fibrewise Wick product $\fibwick$ is compatible
with this $^*$-involution, i.e. 
$(a \fibwick b)^* = (-1)^{\dega a \dega b} b^* \fibwick a^*$ for all 
$a,b \in \WLEnd$. From the unique characterizations of $\rStrwick$ and
$\rwick$ the following lemma is straightforward, see also 
\cite[Lem.~7]{waldmann:2002b}. 
\begin{lemma}
    Let $\Omega =\cc{\Omega}$ be real, then
    \begin{equation}\label{eq:littlerinvolution}
        (\rStrwick)^* = \rStrwick, \qquad \cc{\rwick} = \rwick, \qquad
        (\rEwick)^* = - \rEwick
    \end{equation}
    and hence $(\DStrwick a)^* = \DStrwick a^*$ and $\cc{\Dwick b} =
    \Dwick \cc{b}$ for $a \in \WLEnd$ and $b \in \WL$. Moreover we
    have $(\tauStrwick(A))^* = \tauStrwick (A^*)$ and
    $\cc{\tauwick(f)}= \tauwick(\cc{f})$ and hence $\starStrwick$ and
    $\starwick$ are Hermitian deformations, i.e.
    \begin{equation}\label{eq:HermDef}
        (A \starStrwick B)^* = B^* \starStrwick A^* 
        \quad\textrm{and}\quad \cc{f \starwick g} = \cc{g} \starwick \cc{f}.
    \end{equation}
\end{lemma}
In a next step we extend the fibre metric $h$ to a metric on
$\WLE$ with values in $\WL$ by defining
\begin{equation}
    \label{eq:Hdef}
    H(f\otimes \alpha\otimes s, g \otimes \beta\otimes s')
    :=
    \cc{f}\fibwick g \otimes \cc{\alpha} \wedge \beta h(s,s')
\end{equation}
The following properties are 
immediate\footnote{In \cite{waldmann:2002b} an obvious sign was missing.}:
\begin{lemma}
    For $a \in \WLEnd$, $\Psi,\Psi' \in \WLE$ and $b \in \WL$ we have
    \begin{equation}
        \label{eq:Hcomp1}
        H(a \fibwick \Psi, \Psi') 
        = (-1)^{\dega a \dega \Psi} H(\Psi, a^*\fibwick \Psi'),
    \end{equation}
    \begin{equation}
        \label{eq:Hcomp2}
        H(\Psi, \Psi' \fibwick b) 
        = H(\Psi,\Psi')\fibwick b,
    \end{equation}
    \begin{equation}
        \label{eq:HHStrherm}
        H(\Psi,\Psi')= (-1)^{\dega\Psi\dega\Psi'}\cc{H(\Psi',\Psi)}.
    \end{equation}
\end{lemma}
Now a simple computation using the preceding two lemmas yields the
following compatibility of $\Dwick$ with $H$
\begin{equation}
    \label{eq:DHcomp}
    \Dwick(H(\Psi,\Psi')) 
    = H(\DEwick\Psi, \Psi') + (-1)^{\dega\Psi} H(\Psi,\DEwick \Psi').
\end{equation}
As a consequence we can define a deformed Hermitian metric
$\boldsymbol{h}$ by
\begin{equation}
    \label{eq:hdeformdef}
    \boldsymbol{h}( s, s')
    := \sigma (H(\tauEwick( s), \tauEwick( s'))).
\end{equation}
\begin{proposition}
    \label{proposition:HHStrProper}
    The map $\boldsymbol{h}$ is $\mathbb{C}[[\lambda]]$-linear in the
    second argument and satisfies
    \begin{eqnarray}
        \label{eq:positHerm}
        \boldsymbol{h}( s, s') = \cc{\boldsymbol{h}( s',s)}
        &\textrm{and}&
        \boldsymbol{h}( s, s) \geq 0\\
        \label{eq:bulcomp}
        \boldsymbol{h}( s,  s'\biwick f) 
        & = & \boldsymbol{h}( s,  s')\starwick f\\
        \label{eq:bulStrcomp}
        \boldsymbol{h}(A \biStrwick s, s')
        & = & \boldsymbol{h}( s, A^*\biStrwick s')
    \end{eqnarray}
    for all $ s, s' \in \Gamma^\infty(E)[[\lambda]]$, 
    $f \in C^\infty(M)[[\lambda]]$ and $A \in
    \Gamma^\infty(\End(E))[[\lambda]]$.
\end{proposition}
Here the positivity of $\boldsymbol{h}$ is understood in the sense of
\cite{bursztyn.waldmann:2000b}.

From the properties of $\fibwick$ it is obvious that
$\boldsymbol{h}(s,s')=\sigma(H(\pi_{\cc{z}}\tauEwick(s),\pi_{\cc{z}}\tauEwick(s')))$
and this implies the following lemma using the statements of
Lemma~\ref{lemma:projundef}.
\begin{lemma}
    \label{lemma:metriclocholantihol}
    Let $\mathcal{U}\subseteq M$ be an open subset of $M$.
    In case $s$ or $s'$ are locally holomorphic sections on
    $\mathcal{U}$ we have 
    \begin{equation}
        \label{eq:ssStrlochol}
        \boldsymbol{h}(s,s')|_{\mathcal{U}} =
        h(s,s')|_{\mathcal{U}}.
    \end{equation}
\end{lemma}
To conclude this section we now want to give an explicit local
expression for $\boldsymbol{h}$.
\begin{proposition}
    In case $E$ is an anti-holomorphic vector bundle and $e_\alpha$ is
    an anti-holomorphic frame on $\mathcal{U}_\alpha$ one has
    \begin{equation}
        \label{eq:Metricsonantihol}
        \boldsymbol{h}(s,s')|_{\mathcal{U}_\alpha}
        =\cc{ s^i_\alpha} \starwick \boldsymbol{h}(e_{\alpha,i},
        e_{\alpha,j}) \starwick {s'}^j_\alpha,
    \end{equation}
    where we have written 
    $s|_{\mathcal U_\alpha}= e_{\alpha,j} s^j_\alpha $ and 
    $s'|_{\mathcal U_\alpha}= e_{\alpha,j} {s'}^j_\alpha$.
\end{proposition}
\begin{proof}
    The proof is straightforward using
    Proposition~\ref{proposition:HHStrProper},
    Lemma~\ref{lemma:metriclocholantihol} and the Wick type property
    of $\biwick$.
\end{proof}
\begin{remark}
    In fact it is easy to see that the expression for
    $\boldsymbol{h}(s,s')|_{\mathcal{U}_\alpha}$ given in the
    preceding proposition globally defines the deformed metric
    $\boldsymbol{h}$ by similar arguments as in Section~\ref{sec:local}.
\end{remark}
\begin{proposition}
    In case $E$ is a holomorphic vector bundle and $e_\alpha$ is a
    holomorphic frame on $\mathcal{U}_\alpha$ and $e^\alpha$ is the dual
    frame of the vector bundle $E^*$ on $\mathcal{U}_\alpha$ one has
    \begin{equation}
        \label{eq:Metricsonhol1}
        \boldsymbol{h}(s,s')|_{\mathcal{U}_\alpha}
        = \frac{1}{k^2} 
        h(e_{\alpha,i}, ((s \otimes e^{\alpha,i})^*
        \starStrwick (s' \otimes e^{\alpha,j}))e_{\alpha,j}).
    \end{equation}
\end{proposition}
\begin{proof}
    Again the proof is easily done using
    Proposition~\ref{proposition:HHStrProper},
    Lemma~\ref{lemma:metriclocholantihol} and the Wick type property
    of $\biStrwick$.
\end{proof}
\begin{remark}
    Again the above expression globally defines the metric
    $\boldsymbol{h}$ and here $\boldsymbol{h}$ turns
    out to be canonical in the sense that given a deformation of the
    metric $h$ that satisfies (\ref{eq:ssStrlochol}) and
    (\ref{eq:bulStrcomp}) it coincides with $\boldsymbol{h}$.
\end{remark}

%
%

\section{Morita equivalence of Wick and anti Wick products}
\label{sec:morita}

%
%

In this section we fix a series of closed two-forms $\Omega$ of type
$(1,1)$. Then consider the Wick and anti Wick star product constructed
out of $\Omega$ as in Theorem~\ref{theorem:FedosovI}. Their
characteristic classes satisfy
\begin{equation}
    \label{eq:relclasswaw}
    c(\starwick) - c(\starawick) = -[R^{L_{\mathrm{can}}}] 
    = 2 \pi \im c_1(L_{\mathrm{can}}),
\end{equation}
whence $\starwick$ and $\starawick$ are known to be Morita equivalent
\cite[Thm.~3.1]{bursztyn.waldmann:2002a}. Moreover, as the difference
of their characteristic classes is given by the Chern class of
$L_{\mathrm{can}}$ a Morita equivalence bimodule is obtained by
deforming $\Gamma^\infty(L_{\mathrm{can}})$ into a bimodule for
$\starwick$ from the left and $\starawick$ from the right. Finally,
such a bimodule deformation necessarily exists. However, the concrete
bimodule structure usually depends on non-canonical choices, even
within its equivalence class of bimodule deformations, see e.g. the
constructions in
\cite{bursztyn.waldmann:2002a,bursztyn.waldmann:2000b} as well as the
discussion in \cite{bursztyn.waldmann:2002:pre}.  Thus one may ask the
question whether in our particular situation there is a
\emph{canonical construction}, i.e. only using the Kähler geometry, of
a Morita equivalence bimodule structure for $L_{\mathrm{can}}$. As we
shall show now this is indeed the case.

First we consider the space $\WLL$ where
$\mathcal{L} = \Gamma^\infty(L_{\mathrm{can}})$. We equip $\WL$
with the fibrewise products $\fibwick$ and $\fibawick$ as before. But
for $\WLL$ we use \emph{different} bimodule multiplications, namely
\begin{equation}
    \label{eq:modbimod}
    a \fibbiwick \Psi = S^{-1} a \fibweyl \Psi
    \quad
    \textrm{and}
    \quad
    \Psi \fibbiawick b = \Psi \fibweyl Sb,
\end{equation}
for $a,b \in \WL$ and $\Psi \in \WLL$. Then we indeed obtain a
bimodule as a simple computation shows:
\begin{lemma}
    \label{lemma:modbimod}
    Using $\fibbiwick$ and $\fibbiawick$ the space $\WLL$ becomes a
    bimodule for $(\WL, \fibwick)$ from the left and for
    $(\WL, \fibawick)$ from the right.
\end{lemma}

Now let $\rwick$ and $\rawick$ be the curvature elements as in
Theorem~\ref{theorem:FedosovI}, where we assume to have the
\emph{same} $\Omega$. Moreover, as
connection $\nabla^{L_{\mathrm{can}}}$ we use the canonical connection
induced by the Kähler connection, see (\ref{eq:NablaLcan}). Then we
define $\mathcal{D}^{L}: \WLL \to \WL^{+1}\!\otimes\!\mathcal{L}$
explicitly by
\begin{equation}
    \label{eq:DLDef}
    \mathcal{D}^L \Psi 
    := -\delta \Psi + D^L \Psi + \frac{\im}{\lambda}
    \left(
        \rwick \fibbiwick \Psi 
        - (-1)^{\dega\Psi} \Psi \fibbiawick \rawick
    \right).
\end{equation}
Since the lowest orders of $\rwick$ and $\rawick$ coincide the map
$\mathcal{D}^L$ is well-defined, i.e. it does not produce negative
powers of $\lambda$.
\begin{lemma}
    \label{lemma:DLDerivation}
    The map $\mathcal{D}^L$ is a module derivation along $\Dwick$ and
    $\Dawick$, respectively, i.e.
    \begin{equation}
        \label{eq:DLDerivDwick}
        \mathcal{D}^L(a \fibbiwick \Psi)
        = \Dwick a \fibbiwick \Psi 
        + (-1)^{\dega a} a \fibbiwick \mathcal{D}^L\Psi,
    \end{equation}
    \begin{equation}
        \label{eq:DLDerivDawick}
        \mathcal{D}^L(\Psi \fibbiawick b)
        = \mathcal{D}^L \Psi \fibbiwick b
        + (-1)^{\dega \Psi} \Psi \fibbiwick \Dawick b.
    \end{equation}
\end{lemma}
\begin{proof}
    This is a straightforward computation using the properties of
    $\delta$, $D^L$, and $S$.
\end{proof}
\begin{lemma}
    \label{lemma:DLDL}
    For the symplectic curvature $R$ of the Kähler connection we have
    \begin{equation}
        \label{eq:RPsiPsiR}
        R \fibweyl \Psi = R \fibbiwick \Psi + \lambda \varrho \Psi
        \quad
        \textrm{and}
        \quad
        \Psi \fibweyl R = \Psi \fibbiawick R - \lambda \varrho \Psi,
    \end{equation}
    \begin{equation}
        \label{eq:DLDL}
        (D^L)^2 = \frac{\im}{\lambda} 
        \left(R \fibbiwick \Psi - \Psi \fibbiawick R\right).
    \end{equation}
\end{lemma}
\begin{proof}
    The equation (\ref{eq:RPsiPsiR}) follows directly from the
    definitions and $SR = R + \lambda\varrho$. Moreover,
    (\ref{eq:Ricciform}) together with Lemma~\ref{lemma:Dsquare}
    applied for $E=L_{\mathrm{can}}$ and (\ref{eq:RPsiPsiR}) implies
    (\ref{eq:DLDL}).
\end{proof}
\begin{theorem}
    \label{theorem:Morita}
    We have $(\mathcal{D}^L)^2 = 0$ and
    \begin{equation}
        \label{eq:sigmaLcan}
        \sigma: \ker \mathcal{D}^L \cap \WLO^0\!\otimes\!\mathcal{L}
        \to \Gamma^\infty(L_{\mathrm{can}})[[\lambda]] 
    \end{equation}
    is a $\mathbb{C}[[\lambda]]$-linear bijection with inverse
    denoted by $\tau^{L_{\mathrm{can}}}$. Hence
    \begin{equation}
        \label{eq:leftLcan}
        f \Lbiwick s 
        = \sigma(\tauwick(f) \fibbiwick \tau^{L_{\mathrm{can}}}(s))
    \end{equation}
    \begin{equation}
        \label{eq:rightLcan}
        s \Lbiawick g 
        = \sigma(\tau^{L_{\mathrm{can}}}(s) \fibbiawick \tauawick(g))
    \end{equation}
    defines a bimodule structure on
    $\Gamma^\infty(L_{\mathrm{can}})[[\lambda]]$ for $\starwick$ from
    the left and $\starawick$ from the right, deforming the classical
    bimodule structure of $\Gamma^\infty(L_{\mathrm{can}})$.
\end{theorem}
\begin{proof}
    To compute $(\mathcal{D}^L)^2 = 0$  one has to use
    Lemma~\ref{lemma:DLDL} and the particular properties of $\rwick$
    and $\rawick$ as in Theorem~\ref{theorem:FedosovI}. Taking into
    account that we have used the \emph{same} $\Omega$  for $\rwick$
    and $\rawick$ it easily follows that $(\mathcal{D}^L)^2 = 0$. The
    fact that (\ref{eq:sigmaLcan}) is a bijection is proved in the
    usual fashion by examining the fixed point equation
    \begin{equation}
        \label{eq:tauLfixedpoint}
        \tau^L(s) 
        = s + \delta^{-1}\left(
            D^L\tau^L(s) 
            + \frac{\im}{\lambda} \rwick \fibbiwick \tau^L(s) 
            - \frac{\im}{\lambda} \tau^L(s) \fibbiawick \rawick
        \right)
    \end{equation}
    for $s \in \Gamma^\infty(L_{\mathrm{can}})[[\lambda]]$
    of the strictly contracting operator defined by the right hand
    side, analogously to \cite[Thm.~3]{waldmann:2002b}. Then
    (\ref{eq:leftLcan}) and (\ref{eq:rightLcan}) obviously define a
    bimodule structure deforming the classical one.
\end{proof}
\begin{corollary}
    \label{corollary:Morita}
    The bimodule 
    $(\Gamma^\infty(L_{\mathrm{can}})[[\lambda]], \Lbiwick, \Lbiawick)$
    is a Morita equivalence bimodule for $\starwick$ and $\starawick$.
\end{corollary}
\begin{proof}
    This follows from the general argument in
    \cite[Prop.~1]{waldmann:2002b}.
\end{proof}
\begin{remark}
    \label{remark:Morita}
    Clearly, the above construction is canonical in so far as it uses
    only the Kähler geometry. Hence there exists a distinguished
    bimodule structure on $\Gamma^\infty(L_{\mathrm{can}})[[\lambda]]$
    as desired. Note however, that $\Lbiwick$ and $\Lbiawick$ do not
    have any separation of variables properties.
\end{remark}
\begin{remark}
    \label{remark:WeylMorita}
    Surprisingly, the Weyl product for $\Omega$ is \emph{not} Morita
    equivalent to $\starwick$ or $\starawick$ in general as their
    relative class is given by  
    $c(\starweyl) - c(\starawick) =
    -\frac{1}{2}[R^{L_{\mathrm{can}}}] 
    = \pi\im c_1(L_{\mathrm{can}})$.
    Thus their relative class needs not to be $2\pi\im$ integral in
    general.
\end{remark}

%
%

\appendix

\section{Some Kähler geometry}
\label{sec:kaehler}

In this appendix we shall recall some basic structures on Kähler
manifolds in order to set up our notation and specify our sign
conventions, see e.g.~\cite{wells:1980,kobayashi.nomizu:1969} for
details.

Let $(M, I, g, \omega)$ be a Kähler manifold where $I$ denotes the
complex structure, $g$ the Kähler metric, and $\omega$ the
symplectic Kähler form. Since $I$ is integrable we have holomorphic
local coordinates around any point. If $z^1, \ldots, z^n$ are such
holomorphic local coordinates on $\mathcal{U} \subseteq M$ then we use
the abbreviation $Z_k = \frac{\partial}{\partial z^k}$ and
$\cc{Z}_\ell = \frac{\partial}{\partial \cc{z}^\ell}$. Locally, the
$Z_k$ span the eigenbundle of $I$ to eigenvalue $+\im$ while the
$\cc{Z}_\ell$ span the eigenbundle of $I$ to eigenvalue $-\im$. In
such holomorphic coordinates the Kähler metric and the Kähler form are
given by
\begin{equation}
    \label{eq:localgomega}
    g\big|_{\mathcal{U}} 
    = \frac{1}{2} \, g_{k\cc{\ell}} \, 
    dz^k \vee d\cc{z}^\ell
    \quad
    \textrm{and}
    \quad
    \omega\big|_{\mathcal{U}} 
    = \frac{\im}{2} \, g_{k\cc{\ell}} \, 
    dz^k \wedge d\cc{z}^\ell.
\end{equation}
By $\bigwedge^{(p,q)}T^*M$ we denote the bundle of $p+q$ forms of
type $(p,q)$. In particular,
\begin{equation}
    \label{eq:LcanDef}
    L_{\mathrm{can}} := \mbox{$\bigwedge$}^{(n,0)} T^*M
\end{equation}
is the so-called \emph{canonical line bundle} of 
`holomorphic $n$-forms'. Locally, a section 
$s \in \Gamma^\infty(L_{\mathrm{can}})$ can be written as
$s\big|_{\mathcal{U}} = f dz^1 \wedge \cdots \wedge dz^n$ with  
$f \in C^\infty(\mathcal{U})$. The Kähler connection $\nabla$ extends
to a connection $\nabla^{L_{\mathrm{can}}}$ for $L_{\mathrm{can}}$
which locally is given by
\begin{equation}
    \label{eq:NablaLcan}
    \nabla^{L_{\mathrm{can}}}_X s 
    = (X(f) - dz^k(X) \Gamma^\ell_{k\ell} f) 
    dz^1 \wedge \cdots \wedge dz^n.
\end{equation}
Here $\Gamma^\ell_{km} = dz^\ell(\nabla_{Z_k}Z_m)$ are the Christoffel
symbols of the Kähler connection. As $L_{\mathrm{can}}$ is a line
bundle the curvature of $\nabla^{L_{\mathrm{can}}}$ is just a two-form
$R^{L_{\mathrm{can}}} \in \Gamma^\infty(\bigwedge^{(1,1)} T^*M)$ which
is of type $(1,1)$. Locally one has
\begin{equation}
    \label{eq:RLcan}
    R^{L_{\mathrm{can}}} 
    = \Gamma^j_{kj,\cc{\ell}} \, dz^k \wedge d\cc{z}^\ell
    = - R^j_{jk\cc{\ell}} \, dz^k \wedge d\cc{z}^\ell,
\end{equation}
where $R^j_{mk\cc{\ell}} = dz^j(\hat{R}(Z_k, \cc{Z}_\ell)Z_m)$ are the
components of the curvature tensor $\hat{R}$ of $\nabla$. The
\emph{Ricci form} $\varrho$ is defined by
\begin{equation}
    \label{eq:Ricciform}
    \varrho = \frac{\im}{2} \, R^{L_\mathrm{can}} 
    = - \frac{\im}{2} \, R^j_{jk\cc{\ell}} \, dz^k \wedge d\cc{z}^\ell,
\end{equation}
whence we obtain the following relations for the (first) Chern class
$c_1(L_{\mathrm{can}})$ of the canonical line bundle
\begin{equation}
    \label{eq:ChernclassLcan}
    c_1(L_{\mathrm{can}}) = \frac{\im}{2\pi} [R^{L_\mathrm{can}}]
    = \frac{1}{\pi} [\varrho] \in \HdR^2 (M, \mathbb{Z}).
\end{equation}
Note that $\frac{1}{\pi}\varrho$ is integral but
$\frac{1}{2\pi}\varrho$ does \emph{not} need to be integral.

The symplectic curvature tensor 
$R \in \Gamma^\infty(\bigvee^{(1,1)}T^*M\otimes\bigwedge^{(1,1)}T^*M)$
is defined as usual by
\begin{equation}
    \label{eq:symplRDef}
    R(X,Y,Z,W) = \omega(X, \hat{R}(Z,W)Y)
\end{equation}
for $X,Y,Z,W \in \Gamma^\infty(TM)$ and locally one has
\begin{equation}
    \label{eq:symplRlocal}
    R 
    = \frac{\im}{2} \, g_{k\cc{m}} R^{\cc{m}}_{\cc{\ell} i \cc{j}}
    \, dz^k \vee d\cc{z}^\ell \otimes dz^i \wedge d\cc{z}^j.
\end{equation}
The fact that the Kähler connection is metric implies 
$g_{k\cc{m}} R^{\cc{m}}_{\cc{\ell} i \cc{j}}
= - g_{m\cc{\ell}} R^m_{ki\cc{j}}$
whence one easily obtains
\begin{equation}
    \label{eq:DeltafibRRicci}
    \Deltafib R = \varrho.
\end{equation}
We note that $\cc{R} = R$ and $\cc{\varrho} = \varrho$ are real tensor
fields while $\cc{R^{L_{\mathrm{can}}}} = - R^{L_{\mathrm{can}}}$ is
imaginary.

Now let $E \to M$ be a holomorphic vector bundle of fibre dimension
$k$ with Hermitian fibre metric $h$. By $e = (e_1, \ldots, e_k)$ we
denote a local holomorphic frame of $E$, i.e. 
$e_i \in \Gamma^\infty(E\big|_{\mathcal{U}})$ are holomorphic local
base sections of $E$. Any local section 
$s \in \Gamma^\infty(E\big|_{\mathcal{U}})$ can be 
written as $s = e_i s^i$ with unique local smooth functions 
$s^i \in C^\infty(\mathcal{U})$. Then $s$ is holomorphic if and only
if the $s^i$ are holomorphic. Different local holomorphic frames
$e_\alpha$ and $e_\beta$ on $\mathcal{U}_\alpha$ and
$\mathcal{U}_\beta$, respectively, give rise to holomorphic transition
matrices 
$\phi_{\alpha\beta} = \phi_{\beta\alpha}^{-1}
\in M_k(C^\infty(\mathcal{U}_\alpha\cap\mathcal{U}_\beta))$ by
$e_\beta = e_\alpha \phi_{\alpha\beta}$. As usual they satisfy the
co-cycle identity
\begin{equation}
    \label{eq:cocycle}
    \phi_{\alpha\beta} \phi_{\beta\gamma} \phi_{\gamma\alpha} = \id
\end{equation}
on triple overlaps
$\mathcal{U}_\alpha\cap\mathcal{U}_\beta\cap\mathcal{U}_\gamma$. The
coefficients of a section $s$ transform according to 
$s_\alpha = \phi_{\alpha\beta} s_\beta$.

A connection $\nabla^E$ for $E$ gives rise to a matrix $A$ of local
connection one-forms $A \in M_k (\Gamma^\infty(T^*\mathcal{U}))$ with
respect to a local holomorphic frame $e$ by
\begin{equation}
    \label{eq:ConnectionOneForm}
    \nabla^E_X e = -\im e A(X)
    \quad
    \textrm{for}
    \quad
    X \in \Gamma^\infty(TM).
\end{equation}
Then $\nabla^E$ is called compatible with the holomorphic structure if
$\nabla^E_X s = 0$ for all locally holomorphic sections $s$ and all
vector fields $X$ of type $(0,1)$. Equivalently, the connection
one-forms $A$ with respect to any holomorphic frame are of type
$(1,0)$. The connection is called compatible with the Hermitian fibre
metric if
\begin{equation}
    \label{eq:nablaComph}
    X(h(s,s')) = h(\nabla^E_{\cc{X}} s, s') + h(s, \nabla^E_X s')
\end{equation}
for all $X \in \Gamma^\infty(TM)$ and $s, s' \in
\Gamma^\infty(E)$. This gives the local condition 
$dH = \im(A^*H - HA)$, where $H \in M_k(C^\infty(\mathcal{U}))$ is the
local Hermitian matrix $H_{ij} = h(e_i, e_j)$ defined by the frame
$e$. It is well-known that there exists a unique connection $\nabla^E$
which is compatible with both structures, the holomorphic structure
and the Hermitian fibre metric, see 
e.g.~\cite[Chap.~III, Sec.~2]{wells:1980}. We shall refer to this
connection as the \emph{canonical connection} of $(E, h)$. Locally in
a holomorphic frame one has
\begin{equation}
    \label{eq:LocalCanCon}
    A = \im H^{-1} \partial H.
\end{equation}
Moreover, the curvature tensor $R^E \in \Gamma^\infty(\bigwedge^2 T^*M
\otimes \End(E))$ is of type $(1,1)$. Finally, analogous statements
hold for anti-holomorphic vector bundles as well.

%
%

\begin{footnotesize}

\end{footnotesize}

\end{document}